# Coordinating Multiple Sources for Service Restoration to Enhance Resilience of Distribution Systems

Ying Wang, *Student Member, IEEE*, Yin Xu, *Senior Member, IEEE*, Jinghan He, *Senior Member, IEEE*, Chen-Ching Liu, *Fellow, IEEE*, Kevin P. Schneider, *Senior Member, IEEE*, Mingguo Hong, *Member, IEEE* and Dan T. Ton

*Abstract*—When a major outage occurs on a distribution system due to extreme events, microgrids, distributed generators, and other local resources can be used to restore critical loads and enhance resiliency. This paper proposes a decision-making method to determine the optimal restoration strategy coordinating multiple sources to serve critical loads after blackouts. The critical load restoration problem is solved by a two-stage method with the first stage deciding the post-restoration topology and the second stage determining the set of loads to be restored and the outputs of sources. In the second stage, the problem is formulated as a mixed-integer semidefinite program. The objective is maximizing the number of loads restored, weighted by their priority. The unbalanced three-phase power flow constraint and other operational constraints are considered. An iterative algorithm is proposed to deal with integer variables and can attain the global optimum of the critical load restoration problem by solving a few semidefinite programs in most cases. The effectiveness of the proposed method is validated by numerical simulation with the modified IEEE 13-node test feeder and the modified IEEE 123-node test feeder under plenty of scenarios. The results indicate that the optimal restoration strategy can be determined efficiently in most scenarios.

## I. NOMENCLATURE

*Notation associated with the optimization models:*
Sets and parameters:

| | |
|---|---|
| $\mathcal{N}, \mathcal{L}, \mathcal{G}$ | Sets of all buses, load buses and DG buses in the post-event distribution system, $\mathcal{L}, \mathcal{G} \subseteq \mathcal{N}$ |
| $\mathcal{E}', \mathcal{E}$ | Sets of all available lines and the selected lines to be energized in the target island |
| $\alpha_i$ | Set of phases of bus $i$ |
| $\alpha_{ij}$ | Set of phases of line $(i,j) \in \mathcal{E}$ |
| $\Omega_i$ | Set of adjacent buses of bus $i$ |
| $i \to j$ | Directed line from $i$ to $j$ |
| $\mathbf{Z}_{ij}$ | Impedance matrix of line $(i,j)$ |
| $Y_{ik}^{\phi r}$ | Admittance between phase $\phi$ at bus $i$ and phase $r$ at bus $j$ |
| $w_i$ | Weighting factor of loads at bus $i$, $w_i \geq 0$ |
| $S_{\text{rate},i}$ | Rated power of the DG at bus $i$, $S_{\text{rate},i} = P_{\text{rate},i} + \mathbf{i}Q_{\text{rate},i}$ |
| $s_{\text{load},i}^{\phi}, \mathbf{s}_{\text{load},i}$ | $s_{\text{load},i}^{\phi}$ is complex load power of phase $\phi$ at bus $i$, and $s_{\text{load},i}^{\phi} = p_{\text{load},i}^{\phi} + \mathbf{i}q_{\text{load},i}^{\phi}$, $\mathbf{s}_{\text{load},i} := [s_{\text{load},i}^{\phi}]_{\phi \in \alpha_i}$ |
| $r$ | The root node of a tree in graph theory |

Variables:

| | |
|---|---|
| $V_i^{\phi}, \mathbf{v}_i, \mathbf{V}_i$ | $V_i^{\phi}$ is complex voltage of phase $\phi$ at bus $i$, $\mathbf{v}_i := [V_i^{\phi}]_{\phi \in \alpha_i}$, $\mathbf{V}_i := \mathbf{v}_i \mathbf{v}_i^H$ |
| $\mathbf{v}_i^{\alpha_{ij}}, \mathbf{V}_i^{\alpha_{ij}}$ | $\mathbf{v}_i^{\alpha_{ij}} := [V_i^{\phi}]_{\phi \in \alpha_{ij}}$, $\mathbf{V}_i^{\alpha_{ij}} := \mathbf{v}_i^{\alpha_{ij}} \mathbf{v}_i^{\alpha_{ij}H}$ |
| $I_{ij}^{\phi}, \mathbf{i}_{ij}, \mathbf{I}_{ij}$ | $I_{ij}^{\phi}$ is complex current from bus $i$ to bus $j$ of phase $\phi$, $\mathbf{i}_{ij} := [I_{ij}^{\phi}]_{\phi \in \alpha_{ij}}$, $\mathbf{I}_{ij} := \mathbf{i}_{ij} \mathbf{i}_{ij}^H$ |
| $\mathbf{S}_{ij}$ | $\mathbf{S}_{ij} := \mathbf{v}_i^{\alpha_{ij}} \mathbf{i}_{ij}^H$ |
| $s_i^{\phi}, \mathbf{s}_i$ | $s_i^{\phi}$ is complex power injection of phase $\phi$ at bus $i$, and $s_i^{\phi} = p_i^{\phi} + \mathbf{i}q_i^{\phi}$, $\mathbf{s}_i := [s_i^{\phi}]_{\phi \in \alpha_i}$ |
| $\gamma_i$ | Load status: if load $i$ is restored, $\gamma_i = 1$; otherwise, $\gamma_i = 0$ |
| $a_{ij}$ | Line status: if line $(i,j)$ is selected to be energized, $a_{ij} = 1$; otherwise, $a_{ij} = 0$ |
| $b_{ij}$ | If node $i$ is the parent of node $j$, $b_{ij} = 1$; otherwise, $b_{ij} = 0$ |

*Notations associated with the heuristic for stage 1*

| | |
|---|---|
| $\mathbf{G}$ | The graph corresponding to the target island |
| $\mathbf{V}$ | Set of vertices of graph $\mathbf{G}$, $|\mathbf{V}| = |\mathcal{N}|$ |
| $\mathbf{E}$ | Set of edges of graph $\mathbf{G}$, $|\mathbf{E}| = |\mathcal{E}'|$ |

*Notation associated with the iterative algorithm:*

| | |
|---|---|
| $w^k$ | Weighting factor of the $k$th-level loads |
| $n$ | Number of load levels |
| $P_{\text{load},i}$ | Total real power of load at bus $i$ |
| $\boldsymbol{\gamma}_{\text{CLR-sdp}}^*$ | Solution $\boldsymbol{\gamma}$ of a **CLR-sdp** |
| $\boldsymbol{\gamma}_{\text{CLR}}^*$ | Solution $\boldsymbol{\gamma}$ of the primal **CLR** |
| $\gamma_i^*$ | Element $i$ of $\boldsymbol{\gamma}_{\text{CLR-sdp}}^*$ |
| $W_{\text{CLR-sdp}}^*$ | Optimal value of **CLR-sdp** |
| $W_{\gamma_i=1}$ | Objective value after setting all non-integer $\gamma_i^*$ to 0 |

This work was supported partly by the National Natural Science Foundation of China (51807004), partly by the U.S. Department of Energy (DOE), and partly by the Fundamental Research Funds for the Central Universities(2017YJS185). (Corresponding Author: Yin Xu)

Y. Wang, Y. Xu, and J. He are with the School of Electrical Engineering, Beijing Jiaotong University (BJTU), Beijing 100044, China (e-mails: YingWang1992@bjtu.edu.cn, xuyin@bjtu.edu.cn, jhhe@bjtu.edu.cn).

C.-C. Liu is with the Department of Electrical and Computer Engineering, Virginia Polytechnic Institute and State University (VT), Blacksburg, VA 24061, USA (e-mail: ccliu@vt.edu).

K. P. Schneider is with the Pacific Northwest National Laboratory (PNNL) located at the Battelle Seattle Research Center in Seattle, Washington. (e-mail: kevin.schneider@pnnl.gov).

M. Hong is with the ISO New England, Holyoke, MA 01040, USA (e-mail: mhong@iso-ne.com).

D. T. Ton is with the U.S. Department of Energy (DOE) Office of Electricity Delivery and Energy Reliability (OE), Washington, DC 20585, USA (e-mail: dan.ton@hq.doe.gov).

| | |
|---|---|
| $W_{\text{CLR}}^*$ | Optimal value of the primal **CLR** |
| $\mathcal{L}^k$ | Set of the $k$th-level loads |
| $\mathcal{L}_{\gamma_i=1}$ | Set of loads with $\gamma_i^* = 1$ |
| $\mathcal{L}_{\text{ni}}$ | Set of loads with non-integer value of $\gamma_i^*$ |
| $\mathcal{L}_{\text{c1}}, \mathcal{L}_{\text{c2}}, \mathcal{L}_{\text{c3}}$ | Subsets of $\mathcal{L}_{\text{ni}}$ |
| $\mathcal{L}_{\text{c4}}$ | Set of $k$th-level loads with $\gamma_i^* = 0$ |

## II. Introduction

EXTREME events, including natural disasters, deliberate attacks, and accidents, can cause severe damages to distribution systems, resulting in major outages and significant losses [1]-[3]. Therefore, resilience against extreme events is viewed as an essential feature of smart distribution systems [4]. According to [5], resilience can be enhanced by improving system response and recovery. However, it may be difficult for utility power sources to access interrupted loads after extreme events due to damages of transmission and distribution components and substations. Therefore, this paper is concerned with utilization of multiple local resources, such as distributed generators (DGs) and microgrids (MGs), to restore critical loads on distribution feeders.

Research has been conducted on the utilization of DGs [6]-[11], MGs [12]-[20], and other local sources [21], [22] to restore the critical loads in distribution systems after an extreme event. In [6]-[8], islands energized by DGs are formed to serve load, each of which contains one DG. In [6], the restoration task is formulated as a multi-objective optimization problem and the impact of vehicle-to-grid (V2G) facility of the electric vehicles (EVs) on service restoration is investigated. The multi-agent system approach in [7] is based on expert-system rules to allow the autonomous agents to perform their task. The uncertainties in load demands and renewable generations are considered in [8]. In the authors' prior work, approaches to determining restoration strategies using MGs are proposed [12], [13]. By energizing the paths between MGs and critical loads, several electrical islands are formed, each supported by a MG. The dynamic constraints are considered in [12]. In [15]-[19], the restoration strategies are based on sectionalizing the outage area of a distribution system into several MGs. In [16], a framework is developed to dynamically change the boundaries of MGs for fault isolation and load restoration. A sequential service restoration framework is proposed in [18] to generate a sequence of control actions for switches, DGs, and switchable loads to form multiple isolated MGs. In the aforementioned studies, coordination of multiple sources is not considered.

After an extreme event, by interconnecting and coordinating MGs, DGs, and other local sources, one or more aggregated electrical islands can be formed to restore loads. An example is shown in Fig. 1. By coordinating multiple sources, generation resources can be better allocated and hence restore more critical loads.

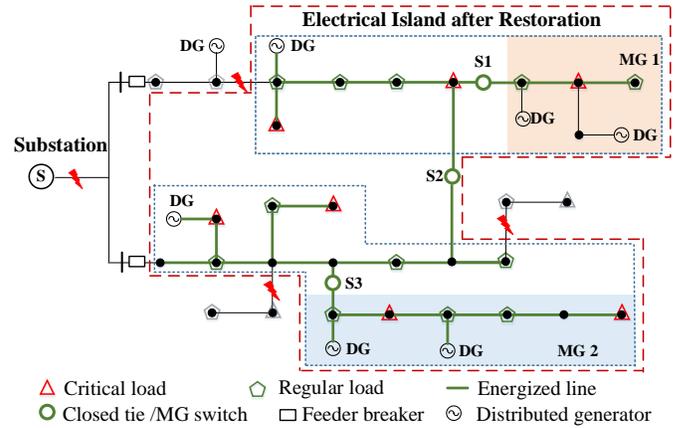

Fig. 1. Coordinating multiple sources for service restoration.

Studies on the coordination of multiple sources for service restoration include the development of planning frameworks, outage management methods, and optimal restoration strategies [23]-[29]. A hierarchical outage management scheme for Multi-Microgrids is developed in [27], which coordinates possible power transfers among the MGs. In [28], the concept of networked MGs for a self-healing distribution system is proposed. Multiple MGs are physically connected via a common bus. In self-healing mode, micro-sources within MGs can be used to support the on-emergency portion of the system. Microgrid formation schemes considering an unbalanced three-phase network model are proposed in [29] to restore critical loads and the coordinated operation of DGs and energy storage (ES) systems is explored. To the best of our knowledge, the existing restoration strategies usually coordinate multiple sources within a MG or at a common bus, while coordinating multiple sources including MGs at multiple locations is rarely considered, which can better explore the benefits of coordination.

In this paper, a decision-making method is proposed to determine the optimal restoration strategy that interconnects and coordinates multiple sources, including DGs and MGs, to serve critical loads. First, the post-restoration topology is decided. Then the set of loads to be restored and the power outputs of DGs are determined. Compared with the state-of-the-art, the major contributions of the paper include:

1) A two-stage critical load restoration method that coordinates multiple sources is proposed, which can better allocate the generation resources and restore more critical loads. By determining the network topology in the first stage and relaxing integer variables representing load status in the second stage, the primal mixed-integer semidefinite program (MISDP) is transformed into a relaxed semidefinite program (SDP) without integer variables, which can be solved efficiently.

2) An iterative algorithm is proposed to deal with non-integer values of load status in the second stage. The algorithm converges in a small number of iterations and the global optimum is attained in most cases. The proposed algorithm significantly improves the computational efficiency compared with the off-the-shelf optimization solvers for mixed-integer programs.

3) A set of sufficient conditions to attain global optimum of the critical load restoration problem are provided and proven: i) the thermal limits of branches are sufficiently large and VAR compensators on critical buses are sufficient; ii) the kW power demands of the loads in the same level are sufficiently different.

The remaining of this paper is organized as follows. Section III introduces the basic formulation of the critical load restoration problem. Section IV presents the two-stage framework for the critical load restoration problem. Section V proposes a mixed-integer semidefinite program for critical load restoration problem and its relaxation to a semidefinite program. Section VI describes the iterative algorithm. Numerical results are presented in Section VII. Section VIII concludes the paper.

### III. PROBLEM FORMULATION

This section formulates the critical load restoration problem. First, the assumptions are presented. Then the idea of coordinating multiple sources for service restoration and the benefits are discussed. Finally, the mathematical formulation of the critical load restoration problem is described.

#### A. Assumptions

It is assumed that after an extreme event, faulted zones are isolated by opening switches or breakers. Survived MGs are disconnected from the main grid and operate in the islanded mode and DGs are disconnected from the distribution network. In addition, suppose that the utility power is unavailable. That is, loads in the distribution system cannot receive power from generators in the transmission systems due to damages of transmission and distribution facilities. As a result, we can only rely on local sources, such as DGs, ESs, and MGs, to restore critical loads.

#### B. Benefits of Coordinating Multiple Sources

In the scenarios described above, one or more electrical islands can be formed by interconnecting MGs, DGs, and loads [30]. These islands can be further connected by closing tie switches between them. For the example shown in Fig. 1, by interconnecting MGs, DGs, and loads, two islands are formed, and then the two islands are further connected by closing tie switch S2 to form one larger island. By doing so, multiple sources, including MGs and different types of DGs, can be coordinated together to realize optimized resources allocation for critical load restoration. Several benefits can be gained by coordinating multiple sources.

1) Multiple power sources at a small capacity can be pooled to serve more loads.

2) Limited generation resources, such as diesel and natural gas, can be optimally allocated to serve critical loads for a longer duration.

3) The control abilities of different types of sources can be fully coordinated to better withstand disturbances during restoration [26].

4) Superior uncertainty management of renewable energy resources can be achieved by coordinating multiple sources with different operational characteristics [23].

#### C. Mathematical Formulation

This paper is focused on exploring the first benefit listed above. A decision-making method is proposed to find the optimal restoration strategy to restore as many critical loads as possible by pooling the capacity of multiple power sources.

Due to the benefits of coordinating multiple sources for service restoration, it is assumed in this paper that outage zones on distribution feeders will be connected to form large islands as long as they are not completely isolated with each other. The islands formed are referred to as the target islands. For each target island, a basic model can be formulated for the critical load restoration problem (**CLR-basic**).

**CLR-basic:**

$$\max f(\boldsymbol{\gamma}) = \sum_{i \in \mathcal{L}} w_i \gamma_i \tag{1}$$

over $\gamma_i \in \{0,1\}$ for $i \in \mathcal{L}$;
$a_{ij}, b_{ij}, b_{ji} \in \{0,1\}$ for $(i,j) \in \mathcal{E}'$;
$V_i^\phi \in \mathbb{C}$, $p_i^\phi, q_i^\phi \in \mathbb{R}^+$ for $i \in \mathcal{N}, \phi \in \alpha_i$;
$I_{ij}^\phi \in \mathbb{C}$ for $(i,j) \in \mathcal{E}, \phi \in \alpha_{ij}$;

s.t.

$$\begin{cases} p_i^\phi + \mathbf{i}q_i^\phi = V_i^\phi \sum_{k \in \Omega_i} I_{ik}^{\phi^H}, \\ I_{ik}^\phi = \sum_{r \in \alpha_k} Y_{ik}^{\phi r}(V_i^\phi - V_k^r) \end{cases} \forall i \in \mathcal{N}, \forall \phi \in \alpha_i \tag{2}$$

$$\begin{cases} 0 \le \sum_{\phi:\phi \in \alpha_i} (p_i^\phi + \gamma_i p_{\text{load},i}^\phi) \le P_{rate,i} \\ 0 \le \sum_{\phi:\phi \in \alpha_i} (q_i^\phi + \gamma_i q_{\text{load},i}^\phi) \le Q_{rate,i} \end{cases}, \forall i \in \mathcal{G} \tag{3}$$

$$\begin{cases} p_i^\phi = -\gamma_i p_{\text{load},i}^\phi \\ q_i^\phi = -\gamma_i q_{\text{load},i}^\phi \end{cases}, \forall i \in \mathcal{N}/\mathcal{G}, \forall \phi \in \alpha_i \tag{4}$$

$$V_{i,\min}^\phi \le |V_i^\phi| \le V_{i,\max}^\phi, \forall i \in \mathcal{N}, \forall \phi \in \alpha_i \tag{5}$$

$$|I_{ij}^\phi| \le I_{ij,\max}^\phi, \forall (i,j) \in \mathcal{E}, \forall \phi \in \alpha_{ij} \tag{6}$$

$$b_{ij} + b_{ji} = a_{ij}, \forall (i,j) \in \mathcal{E}' \tag{7}$$

$$\sum_{j \in \Omega_i} b_{ij} = 1, \forall i \in \mathcal{N}/\text{r} \tag{8}$$

$$b_{rj} = 0, \forall (1,j) \in \mathcal{E}' \tag{9}$$

where $\mathbb{C}$ denotes the set of complex numbers and $\mathbb{R}^+$ the set of nonnegative real numbers.

The objective is to maximize the number of loads restored, weighted by their priority, as indicated by (1). Constraints (2)-(6) are the operational constraints. Constraint (2) represents the unbalanced three-phase power flow, where the first equation indicates the Kirchhoff's current law and second equation is the Ohm's law. Constraints on injection power of a bus are defined by (3) and (4). Specifically, (3) guarantees that the output power of DGs does not exceed their capacity and (4) indicates that the injection power at buses without DG equals the product of negative load power and load status. Variables with subscripts "max" and "min" are the upper and lower limits of the corresponding variables. Constraints (5) and (6) ensure that node voltages and line currents of each phase are within a preset range. Equations (7)-(9) are the radial topological constraints, i.e., spanning tree constraints [31] for the target island. Constraint (7) says that if node $i$ is the parent of node $j$ ($b_{ij} = 1$) or node $j$ is the parent of node $i$ ($b_{ji} = 1$), line $(i, j)$ is selected to be energized ($a_{ij}$). It can be seen that $b_{ij}$ and $b_{ji}$ cannot both equals 1. The parameter r in (8) and (9) represents the root node in the spanning tree. Equation (8) indicates that except the root node, each node has exactly one parent node, and (9) means that the root node has no parent node.

## IV. A Two-Stage Solution Method for CLR-Basic

The challenges in solving the above model lie in the unbalanced power flow equations and large number of integer variables, making the optimization model nonconvex. The rest of this paper will focus on the solution method for the critical load restoration problem **CLR-basic**.

### A. A Two-Stage Solution Method

Integer variables in **CLR-basic** include load status and variables associated with radial topological constraints. To deal with those integer variables, a two-stage solution method is proposed:

Stage 1: Decide the post-restoration topology of the target island. Integer variables associated with topological constraints are decided in this stage.

Stage 2: Determine the loads to be picked up and the power output of all sources in the target island.

### B. A Heuristic to Determine the Topology in Stage 1

In stage 1, a graph-theoretic method is introduced to decide the radial topology for the target island.

The target island can be modeled as an undirected graph $\boldsymbol{G} = (\boldsymbol{V}, \boldsymbol{E})$, where $\boldsymbol{V}$ is the set of vertices and $\boldsymbol{E}$ is the set of edges. Here $\boldsymbol{V}$ and $\boldsymbol{E}$ correspond to the sets $\mathcal{N}$ and $\mathcal{E}'$ in the optimization model, respectively. The length of each edge in $\boldsymbol{E}$ is defined as the electrical distance, i.e., the average self-impedance of all phases of the corresponding line. The distance between a pair of vertices is the total length of edges on the shortest path between the two vertices. The diameter of a graph is the longest distance among all pair of vertices in $\boldsymbol{G}$ [32]. Consider the graph shown in Fig. 2, and assume that all the edges are of the length of 1 unit. Edges (2,4) and (3,4) represent the lines with switches. To maintain a radial structure, vertex 4 can be connected either to vertex 2 or 3 by edge (2,4) or (3,4). Denote the graph with edge (2,4) as $\boldsymbol{G_1}$ and that with (3,4) as $\boldsymbol{G_2}$. The diameters of $\boldsymbol{G_1}$ and $\boldsymbol{G_2}$ are 2 and 3 units, respectively.

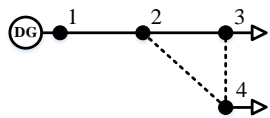

Fig. 2. An example of a graph that represents an electric island.

Some benefits can be gain by choosing a topology with a short diameter. First, the power losses can be small because power does not need to travel over long distance. Second, the voltage profile can be better because reactive sources can be close to reactive loads. Third, the risk of overloading is reduced.

For the example shown in Fig. 2, it can be proven that compared with $\boldsymbol{G_2}$, by choosing $\boldsymbol{G_1}$ as the post-restoration topology, the power loss on line (2,3) is smaller and the voltage drop at nodes 3 and 4 is reduced.

Therefore, the tree with the minimal diameter is chosen as the post-restoration topology. This problem is identical to the minimum diameter spanning tree (MDST) problem in graph theory. The algorithm to find the MDST in [32] is used in this work. The MDST can be found in $\mathcal{O}(mn + n^2 \log n)$ time, where $n = |\boldsymbol{V}|$ and $m = |\boldsymbol{E}|$.

### C. Solution Method for Stage 2

After the post-restoration topology is decided, the topological constraints (7)-(9) can be removed from **CLR-basic**. The set of loads to be restored and the power output of each source are then be determined in this stage. The challenges lie in the nonconvex power flow equation constraints and the integer variables of load status. To overcome these difficulties, the **CLR-basic** is first reformulated as a mixed-integer semidefinite program **CLR** and then relaxed to a semidefinite program **CLR-sdp**, which will be introduced in Section V. Finally, an iterative algorithm is proposed in Section VI to solve **CLR**.

## V. A Relaxed Semidefinite Program Formulation

Note that the critical load restoration problem **CLR-basic** is actually includes some common decision variables and constraints with the optimal power flow (OPF) problem. In particular, both problems adjust the output of DGs to achieve their objectives and consider the power flow equations and operational limits on bus voltages, line currents, and DG power as the constraints. The critical load restoration problem introduces additional topological constraints and load status variables, which make it a mixed-integer nonlinear program. For distribution systems, three relaxed models for the power flow constraints have been proposed, i.e., the second-order cone model for balanced networks [33], [34], the semidefinite model and linear approximation model for unbalanced networks [35]. Any of these models can be used to formulate the critical load restoration problem, as discussed in Appendix A. In this paper, the semidefinite model is chosen, because it is capable of modeling unbalance in three phases, power losses and voltage profile with high accuracy.

By applying the variable transformation proposed in [35], i.e., use the slack variables $\boldsymbol{I}_{ij} = \boldsymbol{i}_{ij}\boldsymbol{i}_{ij}^H$ and $\boldsymbol{V}_i = \boldsymbol{v}_i\boldsymbol{v}_i^H$ to replace the variables $\boldsymbol{i}_{ij}$ and $\boldsymbol{v}_i$, the **CLR-basic** modeled by (1)-(6) is transformed into a mixed-integer semidefinite program, which is referred to as **CLR**.

**CLR:**
$$\max \quad f(\boldsymbol{\gamma})$$
over $\gamma_i \in \{0,1\}$ for $i \in \mathcal{L}$;
$\boldsymbol{s}_i \in \mathbb{C}^{|\alpha_i|}$, $\boldsymbol{V}_i \in \mathbb{H}^{|\alpha_i| \times |\alpha_i|}$ for $i \in \mathcal{N}$;
$\boldsymbol{S}_{ij} \in \mathbb{C}^{|\alpha_{ij}| \times |\alpha_{ij}|}$, $\boldsymbol{I}_{ij} \in \mathbb{H}^{|\alpha_{ij}| \times |\alpha_{ij}|}$ for $i \to j$
s.t. (5)

$$\sum_{k:k\to i} \text{diag}(S_{ki} - Z_{ki}I_{ki}) + s_i = \sum_{j:i\to j} \text{diag}(S_{ij})^{\alpha_i}, \forall i \in \mathcal{N} \quad (10)$$

$$0 \leq \text{sum}(s_i + \gamma_i s_{\text{load},i}) \leq S_{\text{rate},i}, \forall i \in \mathcal{G} \quad (11)$$

$$s_i = -\gamma_i s_{\text{load},i}, \forall i \in \mathcal{N}/\mathcal{G} \quad (12)$$

$$v_{i,\min} \leq \text{diag}(V_i) \leq v_{i,\max}, \forall i \in \mathcal{N} \quad (13)$$

$$\text{diag}(I_{ij}) \leq i_{ij,\max}, \forall (i,j) \in \mathcal{E} \quad (14)$$

$$V_j = V_i^{\alpha_{ij}} - (S_{ij}Z_{ij}^H + Z_{ij}S_{ij}^H) + Z_{ij}I_{ij}Z_{ij}^H, \forall i \to j \quad (15)$$

$$\begin{bmatrix} V_i^{\alpha_{ij}} & S_{ij} \\ S_{ij}^H & I_{ij} \end{bmatrix} \succcurlyeq 0, \forall i \to j \quad (16)$$

$$\text{rank}\begin{bmatrix} V_i^{\alpha_{ij}} & S_{ij} \\ S_{ij}^H & I_{ij} \end{bmatrix} = 1, \forall i \to j \quad (17)$$

where the slack variable $S_{ij} = v_i^{\alpha_{ij}} i_{ij}^H$, $|\cdot|$ denotes the number of elements in a set, $\mathbb{C}^m$ the set of $m$-dimensional complex vectors, $\mathbb{H}^{m \times m}$ the set of $m \times m$ complex Hermitian matrices, and $\mathbb{C}^{m \times m}$ the set of $m \times m$ complex matrices.

Constraint (10) represents the power balance equation, where diag(·) returns a column vector of the main diagonal elements of a matrix. Constraints (11) and (12) correspond to (3) and (4), where sum(·) calculates the sum of all elements in a vector. Constraints (13) and (14) correspond to (5) and (6), respectively, where $v_{i,\min} = [(V_{i,\min}^\phi)^2]_{\phi \in \alpha_i}$ and $v_{i,\max} = [(V_{i,\max}^\phi)^2]_{\phi \in \alpha_i}$. For the reference bus, $V_{i,\min}^\phi = V_{i,\max}^\phi = 1$, while $V_{i,\min}^\phi = 0.95$, $V_{i,\max}^\phi = 1.05$ for other buses. Constraint (15) represents the Ohm's law. Constraints (16) and (17) indicate that the matrix $\begin{bmatrix} V_i^{\alpha_{ij}} & S_{ij} \\ S_{ij}^H & I_{ij} \end{bmatrix}$ is positive semidefinite and has a rank of 1 because the matrix is a product of a non-zero vector and its conjugate transpose, i.e.,

$$\begin{bmatrix} V_i^{\alpha_{ij}} & S_{ij} \\ S_{ij}^H & I_{ij} \end{bmatrix} = \begin{bmatrix} v_i^{\alpha_{ij}} \\ i_{ij}^H \end{bmatrix} \begin{bmatrix} v_i^{\alpha_{ij}} \\ i_{ij}^H \end{bmatrix}^H. \quad (18)$$

It should be noted that the power flow equation constraint (2) is corresponding to (10) and (15)-(17), where (10) and (15) are equivalent to the first and second equations in (2), respectively, and (16)-(17) are derived from the definition of slack variable $S_{ij}$, i.e., $S_{ij} = v_i^{\alpha_{ij}} i_{ij}^H$.

By replacing $v_i v_i^H$ and $i_{ij} i_{ij}^H$ by $V_i$ and $I_{ij}$, the phase angles of bus voltages and line currents are eliminated, that is, the phase angles are relaxed. However, the phase angles can be recovered after the solution is obtained [35].

In order to develop a convex model, two relaxations are applied: First, the binary variables $\gamma_i \in \{0,1\}$ are relaxed to continuous variables $\gamma_i \in [0,1]$; Second, the rank constraint (17) is removed. A semidefinite relaxation of **CLR** is then obtained as follows.

**CLR-sdp:**

$$\max \quad f(\gamma)$$
$$\text{over} \quad \gamma_i \in [0,1] \text{ for } i \in \mathcal{L}, s, V, S, \text{ and } I$$
$$\text{s.t.} \quad (10)-(16)$$

The semidefinite relaxation **CLR-sdp** is exact if its solution satisfies (17) and does not contain any non-integer element in $\gamma$. The results in [35] indicate that the solution of **CLR-sdp** satisfies (17) in most cases. However, the solution $\gamma$ of **CLR-sdp** may contain non-integer elements, which implies that the corresponding loads are partially restored, making the solution infeasible for the primal **CLR**. There are two reasons that may lead to non-integer values in $\gamma$:

1) Insufficient generation capacity, i.e., restoring more load will cause a violation of constraint (11);
2) The node voltages or line currents hit their lower or upper limits, i.e., constraints (13) or (14) will be violated if more load is restored.

It is a challenging task to deal with the non-integer values in $\gamma$ efficiently and find the globally optimal solution to the primal **CLR**. The next section will address this issue.

## VI. AN ITERATIVE ALGORITHM

In the solution of **CLR-sdp**, the load status $\gamma_i$ may be a value between 0 and 1, indicating that the corresponding load is partially restored. However, the load connected to the feeder through one switch can only be restored as a whole. To deal with non-integer values in $\gamma$, an iterative algorithm is proposed in this paper. It can be proved that under two conditions, the solution provided by the proposed algorithm is globally optimal for the primal **CLR**.

### A. Selection of Weighting Factors

In this paper, the loads are divided into $n$ levels according to their priority. Loads in the same level are assigned the same value for their weighting factor. Denote the weighting factor of the $k$th-level loads by $w^k \geq 0$, $k = 1, \ldots, n$. Assume that the loads in the 1st level are most critical while those in the $n$th level least important. Therefore, $w^1 > \cdots > w^n$. For instance, for load $i$ whose level is $k$, its weighting factor $w_i$ equals $w^k$. Two virtual levels are introduced, i.e., the 0th and $(n+1)$th levels, with $w^0 = \infty$ and $w^{n+1} = 0$.

Noted that the number of load levels, the level of each load, and weighting factors of different levels are user defined. To ensure that loads with higher priority are restored before those with lower priority, the weighting factors for loads in different levels should be sufficiently different. That is, the following two conditions should be satisfied.

1) For two levels $k_1 < k_2$, $w^{k_1}/P_{\text{load},i} \gg w^{k_2}/P_{\text{load},j}$ for any load $i \in \mathcal{L}^{k_1}$ and $j \in \mathcal{L}^{k_2}$.
2) For an arbitrary level $j$, $w^j > \sum_{k:k>j} w^k |\mathcal{L}^k|$.

### B. Iterative Algorithm

The pseudo code of the proposed algorithm is shown below. It solves the primal **CLR** by solving **CLR-sdp** iteratively. It returns the values of all variables. If the solution $\gamma$ of **CLR-sdp**, denoted by $\gamma^*_{\text{CLR-sdp}}$, contains non-integer elements, additional constraints on $\gamma$ are added to **CLR-sdp** by the function ADDCONSTRAINTS. The value of $\gamma^*_{\text{CLR-sdp}}$ is then updated by solving **CLR-sdp** with the new constraints. The iteration keeps running until there is no non-integer elements in $\gamma^*_{\text{CLR-sdp}}$.

| Iterative Algorithm for the Primal CLR |
|---|
| 1 Relax the primal **CLR** to **CLR-sdp** |
| 2 Solve **CLR-sdp** |
| 3 $\gamma^*_{\text{CLR-sdp}} \leftarrow$ solution $\gamma$ of **CLR-sdp** |
| 4 *while* $\gamma^*_{\text{CLR-sdp}}$ contains non-integer elements |
| 5     ADDCONSTRAINTS |
| 6     Solve **CLR-sdp** with new constraints |
| 7     $\gamma^*_{\text{CLR-sdp}} \leftarrow$ solution $\gamma$ of the updated **CLR-sdp** |
| 8 *end while* |
| 9 Conduct BFM-OPF in [35] based on $\gamma^*_{\text{CLR-sdp}}$ to obtain the solution of $s, V, S, I$ |
| 10 *return* the solution $\gamma, s, V, S, I$ |

Note that when the condition "$\gamma^*_{\text{CLR-sdp}}$ contains no non-integer elements" is satisfied, the value of $\gamma$ (load status) is actually determined by the additional constraints added to **CLR-sdp**. As a result, in the last iteration of the *while* loop defined in steps 4-8, the objective value of **CLR-sdp** is fixed regardless the output of DGs. Therefore, an OPF is conducted in step 9 to determine the output of DGs with the status of loads indicated by $\gamma^*_{\text{CLR-sdp}}$.

Denote the optimal value of **CLR** and **CLR-sdp** by $W^*_{\text{CLR}}$ and $W^*_{\text{CLR-sdp}}$, respectively. Let $W_{\gamma_i=1}$ denote the value of the objective function $f(\gamma)$ after the non-integer elements in $\gamma^*_{\text{CLR-sdp}}$ are set to 0. Evidently, $W^*_{\text{CLR-sdp}} \geq W^*_{\text{CLR}} \geq W_{\gamma_i=1}$. That is, $W^*_{\text{CLR-sdp}}$ is an upper bound of $W^*_{\text{CLR}}$ while $W_{\gamma_i=1}$ a lower bound. Consider $f(\gamma) = W_{\gamma_i=1}$ as the starting point. Then $W^*_{\text{CLR-sdp}} - W_{\gamma_i=1}$ defines an upper bound on the improvement that can be made in the objective value. Based on this fact, the function ADDCONSTRAINTS first secures loads that are fully restored according to the solution $\gamma^*_{\text{CLR-sdp}}$ of **CLR-sdp**, i.e., adding the following constraints to **CLR-sdp**.

$$\gamma_i = 1 \text{ for } i \in \mathcal{L}_{\gamma_i=1} = \{i | \gamma^*_{\text{CLR-sdp}}(i) = 1\}$$

Then, it determines loads that cannot be restored and restricts the corresponding elements in $\gamma$ to 0. The procedure is as follows.

Step 1: Find the set of loads associated with non-integer elements in $\gamma^*_{\text{CLR-sdp}}$, i.e., $\mathcal{L}_{\text{ni}} = \{i | \gamma^*_{\text{CLR-sdp}}(i) \in (0,1)\}$.

Step 2: Identify the load level $K^* \in \{1, ..., n+1\}$ such that $w^{K^*} \leq W^*_{\text{CLR-sdp}} - W_{\gamma_i=1} < w^{K^*-1}$.

Step 3: Find the subset $\mathcal{L}_{c1}$ of $\mathcal{L}_{\text{ni}}$ such that the weighting factor of the loads in $\mathcal{L}_{c1}$ is equal to or greater than $w^{K^*-1}$. If $\mathcal{L}_{c1} \neq \emptyset$, add the following constraints to **CLR-sdp** and *return*.

$$\gamma_i = 0 \text{ for } i \in \mathcal{L}_{c1}$$

Step 4: Find the subset $\mathcal{L}_{c2}$ of $\mathcal{L}_{\text{ni}}$ such that weighting factor of the loads in $\mathcal{L}_{c2}$ is equal to $w^{K^*}$ and the solved load bus voltage matrix $V_i$ or connected line current matrix $I_{ij}$ for load $i \in \mathcal{L}_{c2}$ satisfies any of the bounds in (13) or (14). If $\mathcal{L}_{c2} \neq \emptyset$, add the following constraints to **CLR-sdp** and *return*.

$$\gamma_j = 0 \text{ for } j = \underset{i}{\operatorname{argmin}} \gamma^*_{i:i \in \mathcal{L}_{c2}}$$

where the function $\underset{x}{\operatorname{argmin}} f(x)$ returns the points $x$ for which $f(x)$ attains its smallest value. That is, load $j$ is the load in $\mathcal{L}_{c2}$ whose $\gamma^*_i$ is the smallest.

Step 5: If $\mathcal{L}_{c2} = \emptyset$, loads in $\mathcal{L}_{\text{ni}}$ must be in the same level $K^*$. Determine the maximum number of loads in $\mathcal{L}_{\text{ni}}$ that can be restored by equation (19) as follows

$$n_{\text{re}} = \text{floor}\left(\frac{W^*_{\text{CLR-sdp}} - W_{\gamma_i=1}}{w^{K^*}}\right) \quad (19)$$

where floor($\cdot$) rounds the input number to the nearest integer smaller than or equal to it. Find the subset $\mathcal{L}_{c3}$ of $\mathcal{L}_{\text{ni}}$ such that $|\mathcal{L}_{c3}| = |\mathcal{L}_{\text{ni}}| - n_{\text{re}}$ and for any $i \in \mathcal{L}_{c3}$ and $j \in \mathcal{L}_{\text{ni}} - \mathcal{L}_{c3}$, $\gamma_i \leq \gamma_j$. Find a set $\mathcal{L}_{c4}$ such that weighting factor of the loads in $\mathcal{L}_{c4}$ is equal to $w^{K^*}$ and $\gamma^*_i = 0$. Add the following constraints to **CLR-sdp** and *return*.

$$\gamma_i = 0 \text{ for } i \in \mathcal{L}_{c3} \cup \mathcal{L}_{c4}$$

Constraints that stated in steps 3-5 deal with three different conditions:

If $\mathcal{L}_{c1} \neq \emptyset$, the constraints mentioned in step 3 will be added, i.e., the load status for all loads $i \in \mathcal{L}_{c1}$ is set to zero. The reasons are as follows. If any load in $\mathcal{L}_{c1}$ is restored, the inequation $f(\gamma) \geq W_{\gamma_i=1} + w^{K^*-1} > W^*_{\text{CLR-sdp}}$ holds, contradicting the fact that $f(\gamma) \leq W^*_{\text{CLR}} \leq W^*_{\text{CLR-sdp}}$. Therefore, no load in $\mathcal{L}_{c1}$ can be restored.

If $\mathcal{L}_{c1} = \emptyset$ and $\mathcal{L}_{c2} \neq \emptyset$, the constraints stated in step 4 will be added. This is because if all loads in $\mathcal{L}_{c2}$, together with loads whose $\gamma^*_{\text{CLR-sdp}}(i) = 1$, are restored, constraints (13) or (14) will be violated. Therefore, at least one load in $\mathcal{L}_{c2}$ cannot be restored. The one with the minimum value of $\gamma_i$ is abandoned.

However, if $\mathcal{L}_{c1} = \emptyset$ and $\mathcal{L}_{c2} = \emptyset$, the algorithm goes to step 5. In this condition, loads in $\mathcal{L}_{\text{ni}}$ cannot be fully restored due to insufficient power capacity of DGs. Since it is assumed that the difference in weighting factors of different load levels are sufficiently large (see Section VI. A), loads of low priority will not be restored until high-priority loads are fully restored. Therefore, loads in $\mathcal{L}_{\text{ni}}$ are in the same level $K^*$. In addition, $K^*$ satisfies $w^{K^*} \leq W^*_{\text{CLR-sdp}} - W_{\gamma_i=1} < w^{K^*-1}$, so $n_{\text{re}} \geq 1$, indicating that some loads in $\mathcal{L}_{\text{ni}}$ may be fully restored if others were not restored. The maximum number of loads that can be fully restored are then determined and the loads with small values of $\gamma_i$ are abandoned.

Flowchart of function ADDCONSTRAINTS is shown in Fig. 3.

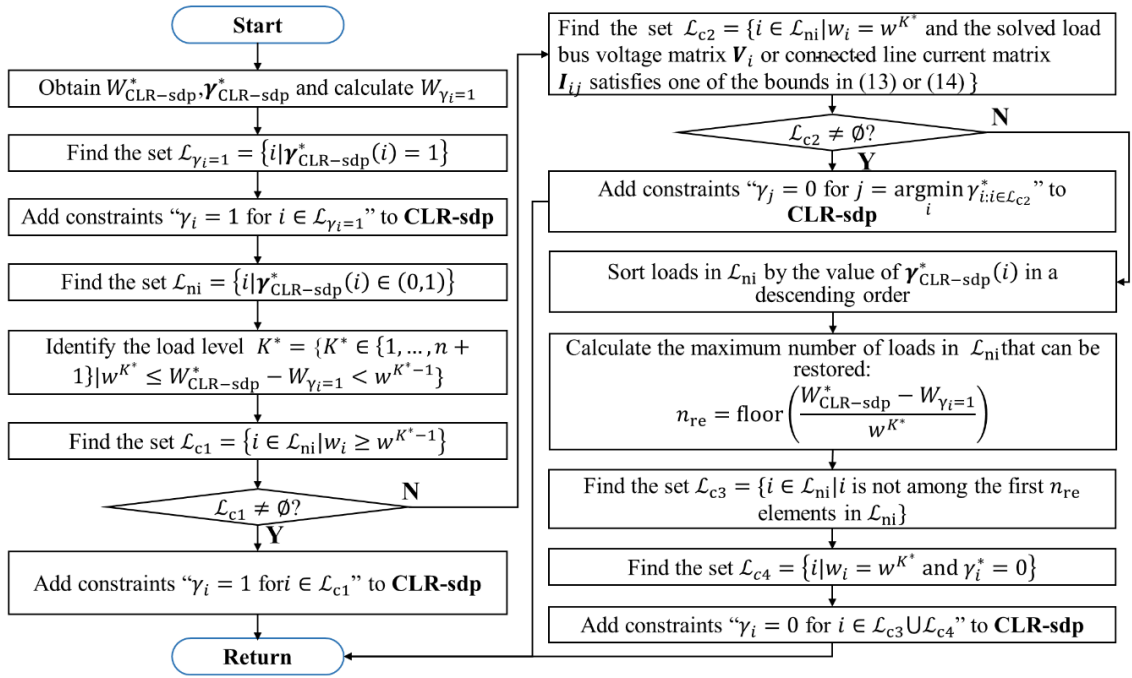

Fig. 3. The flow chart of the function ADDCONSTRAINTS

### C. Convergence Analysis

The iterative process will end when all elements in the solution $\gamma^*_{\text{CLR-sdp}}$ are integers. In the proposed algorithm, the function ADDCONSTRAINTS sets at least one element of $\gamma^*_{\text{CLR-sdp}}$ to a fixed integer value, i.e., 0 or 1, in each iteration. Since the number of elements in $\gamma^*_{\text{CLR-sdp}}$ is finite, the proposed algorithm can determine their values in finite iterations. Therefore, convergence of the proposed algorithm is guaranteed.

### D. Global Optimum Verification

According to the iterative algorithm, the global optimum can be verified based on the solutions of **CLR-sdp** in each iteration. The *optimality criterion* is as follow:

*Optimality criterion*: If no constraint is added in step 4 of the function ADDCONSTRAINTS and $W^{[j+1]}_{\gamma_i=1} \geq W^{[j]}_{\gamma_i=1} + w^{K^*[j]} n^{[j]}_{\text{re}}$ holds for the $j$th iteration, where $j = 0,1, \ldots m-1$, the global optimum of **CLR** is attained.

In the above statement, subscript $[j]$ indicates that the corresponding variable that takes the value obtained in the $j$th iteration. The lines 1-3 of the iterative algorithm are considered as the 0th iteration. The total number of iterations is $m$.

A proof of this criterion is given in Appendix B.

### E. Sufficient Conditions to Attain Global Optimum

According to the *optimality criterion*, the global optimum of **CLR** can be attained if the following conditions hold:

*Condition 1*: the thermal limits of branches are sufficiently large and VAR compensators on critical buses are sufficient.

*Condition 2*: the kW power demands of the loads in the same level are sufficiently different.

If condition 1 holds, the current and voltage constraints will not be activated, so the ADDCONSTRAINTS function will add constraints in step 4. If condition 2 holds, the values of $w_i/P_{\text{load},i}$ for loads in a priority level are sufficiently different, so no constraint is added in step 5 of ADDCONSTRAINTS.

## VII. CASE STUDIES

The proposed restoration model and algorithm have been programmed in MATLAB R2016a with CVX package [36]. The semidefinite programs in the iterative algorithm are optimized using optimizer in MOSEK [37]. The numerical experiments are carried out on a personal computer with Intel Core I5 CPU at 2.40 GHz with 8 GB of RAM.

### A. Case I: The Modified IEEE 13-Node Test Feeder

The IEEE 13-node test feeder [38] is modified to include 2 MGs, as shown in Fig. 4. Two switches 645-632 and 671-692 are added as MG switches. Two tie lines 634-675 and 633-671 are added to provide flexibility in topology. Three DGs and an ES are available and connected at buses 646, 633, 675, and 645, respectively. The detailed information of the DGs and the ES is listed in Table I. Loads are divided into three levels and the weighting factors are 100, 10, and 0.2, respectively. Therefore, $n = 3$, $w^1 = 100$, $w^2 = 10$, and $w^3 = 0.2$. Load levels are marked in Fig.4.

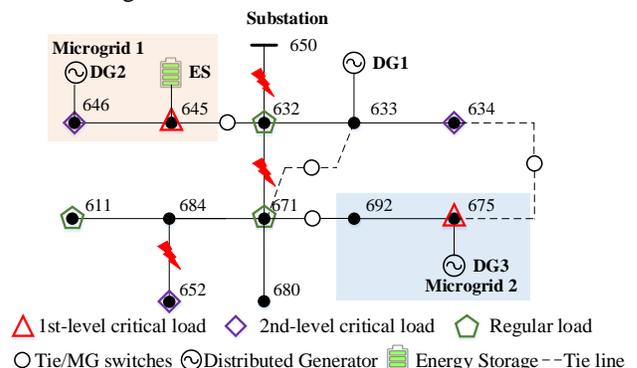

Fig. 4. One-line diagram of modified IEEE 13-node test feeder with 2 MGs.

TABLE I
DGS AND ES DATA

| Generator (#) | Microgrid (#) | Node (#) | Generator Type | Rating (kW) |
|---|---|---|---|---|
| DG1 | / | 633 | Diesel | 600 |
| DG2 | 1 | 646 | Diesel | 200 |
| ES | 1 | 645 | Storage | 280 |
| DG3 | 2 | 675 | Diesel | 360 |

Suppose that two faults occurred and have been isolated properly. The faulted lines are 632-671 and 684-652. Both MGs and DG1 are available for restoration. There are four switches in the system, i.e., 632-645, 633-671, 671-692, and 634-675. The target island $G$ is shown in the dotted box in Fig.5.

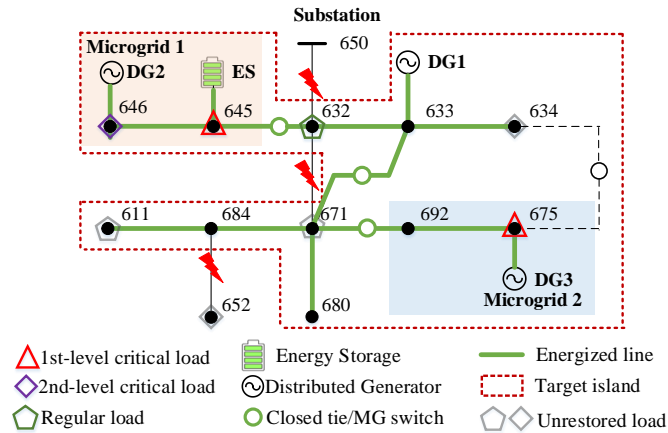

Fig.5. The topology of post-restoration system.

The proposed method is applied to determine an optimal restoration strategy. In the first stage, the MDST represented by green line in Fig. 5 is found. The switches 632-645, 633-671, and 671-692 are closed.

In stage 2, the restoration strategy is obtained in two iterations. The results are listed in Table II.

TABLE II
RESULTS AFTER EACH ITERATION

| Node (#) | $w_i$ | Initial CLR-sdp $\gamma_i^*$ | 1st iteration Constraints | 1st iteration $\gamma_i^*$ | 2nd iteration Constraints | 2nd iteration $\gamma_i^*$ |
|---|---|---|---|---|---|---|
| 632 | 0.2 | 0 | | 1 | **1** | 1 |
| 634 | 10 | 0.49 | **0** | 0 | 0 | 0 |
| 671 | 0.2 | 0 | | 0 | **0** | 0 |
| 675 | 100 | 1 | **1** | 1 | 1 | 1 |
| 645 | 100 | 1 | **1** | 1 | 1 | 1 |
| 646 | 10 | 1 | **1** | 1 | 1 | 1 |
| 611 | 0.2 | 0 | | 0.56 | 0 | 0 |
| $W^*_{\text{CLR-sdp}}$ | | 214.86 | | 210.31 | | 210.2 |
| $W_{\gamma_i=1}$ | | 210 | | 210.2 | | 210.2 |
| ratio | | 3.53e-7 | | 8.96e-6 | | 8.00e-9 |

Take the first iteration as an example to demonstrate the algorithm. In the solution $\gamma^*_{\text{CLR-sdp}}$ of the initial **CLR-sdp**, one non-integer value is found (bus 634). In the 1st iteration, the function ADDCONSTRAINTS is first executed:
1) The set $\mathcal{L}_{\gamma_i=1} = \{675, 645, 646\}$, so the constraints $\gamma_i = 1$ are added for loads at buses 675, 645, and 646.
2) Step 1: The set $\mathcal{L}_{\text{ni}} = \{634\}$.
3) Step 2: The load level $K^* = 3$ is identified by finding $w^{K^*}$ that satisfies $w^{K^*} \leq W^*_{\text{CLR-sdp}} - W_{\gamma_i=1} < w^{K^*-1}$, i.e., $0.2 \leq 4.86 < 10$.
4) Step 3: The set $\mathcal{L}_{c_1} = \{634\}$, and the constraint $\gamma_i = 0$ is added for the load at bus 634.
5) The function ADDCONSTRAINTS returns the added constraint to the **CLR-sdp**.

The entry "ratio" in Table II qualifies how close a **CLR-sdp** is to rank one. To qualify how close a matrix to rank one, obtain the ratio of the two largest eigenvalues $\lambda_1, \lambda_2$ ($|\lambda_1| \geq |\lambda_2|$) of the matrix, i.e., $|\lambda_2|/|\lambda_1|$. The smaller the ratio is, the closer the matrix is to rank one [35]. The maximum ratio over all matrices in (17) is the entry "ratio" in Table II. It can be seen that the solutions of each iteration are numerically exact, i.e., the constraint (17) is satisfied for all solutions.

The results indicate that loads at buses 675, 645, 646, and 632 are restored, as shown in Fig. 5. In addition, the optimal operating state of the post-restoration system is determined. The power outputs of DG1, DG2, DG3, and ES are 552.54kW, 200.00kW, 360.00kW, and 231.76kW, respectively. The power loss is 1.30kW. In this case, the voltage of bus connected with DG1, i.e., bus 633, is selected as the reference with its magnitude equals 1 p.u. and phase angle 0 degree. The three-phase voltage magnitude and angles are listed in Table III.

TABLE III
THREE-PHASE VOLTAGE MAGNITUDE AND ANGLES

| Node (#) | $|V_i^a|$(p.u.) | $\angle V_i^a(°)$ | $|V_i^b|$(p.u.) | $\angle V_i^b(°)$ | $|V_i^c|$(p.u.) | $\angle V_i^c(°)$ |
|---|---|---|---|---|---|---|
| 632 | 1.0006 | -0.0160 | 0.9983 | 120.0107 | 0.9980 | -119.9423 |
| 633 | 1.0000 | 0.0000 | 1.0000 | 120.0000 | 1.0000 | -120.0000 |
| 634 | 1.0000 | 0.0000 | 1.0000 | 120.0000 | 1.0000 | -120.0000 |
| 645 | 1.0006 | -0.0160 | 0.9973 | 120.0995 | 0.9972 | -119.8462 |
| 646 | 1.0006 | -0.0160 | 0.9971 | 120.0938 | 0.9970 | -119.8505 |
| 671 | 0.9994 | -0.0436 | 0.9994 | 119.9592 | 0.9994 | -120.0447 |
| 692 | 0.9994 | -0.0436 | 0.9994 | 119.9592 | 0.9994 | -120.0447 |
| 675 | 0.9976 | -0.0724 | 0.9976 | 119.9265 | 0.9976 | -120.0697 |
| 680 | 0.9994 | -0.0436 | 0.9994 | 119.9592 | 0.9994 | -120.0447 |
| 684 | 0.9994 | -0.0425 | 0.9994 | 119.9562 | 0.9994 | -120.0447 |

According to the *optimality criterion*, the global optimality is attained. Since the number of loads is small, an exhausted search is conducted. The global optimal solution is found, i.e., $\gamma^*_{\text{CLR}} = \{1,0,0,1,1,1,0\}$ and $W^*_{\text{CLR}} = 210.2$, which are the same as the results provided by the proposed method.

B. *Case II: The Modified 123-Node Test System*

1) *System information*

A modified IEEE 123-node test system with 3 MGs [26] is used to test the effectiveness of the proposed method, as shown in Fig. 6. In the normal condition, the 3 MGs operate in the grid-connected mode. After an extreme event, 3 MGs change to the islanded mode. The information of the available DGs in the MGs, including diesel generators, ES, and PVs, is listed in Table IV. Assume that the active power of PVs is constant. Loads are divided into three levels with weighting factors 100, 10, and 0.1, respectively.

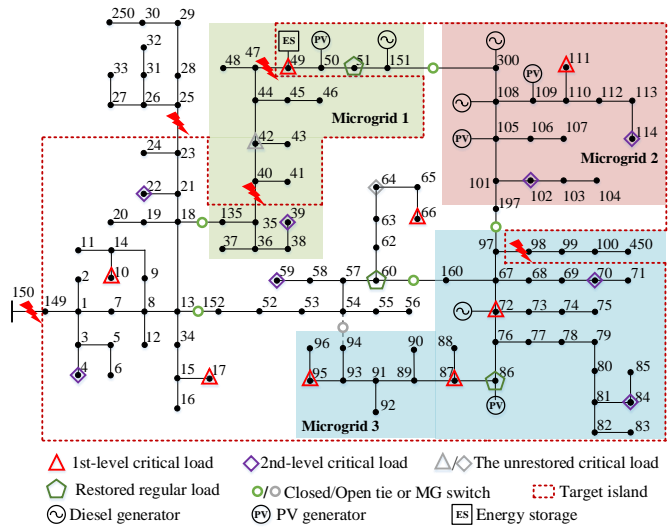

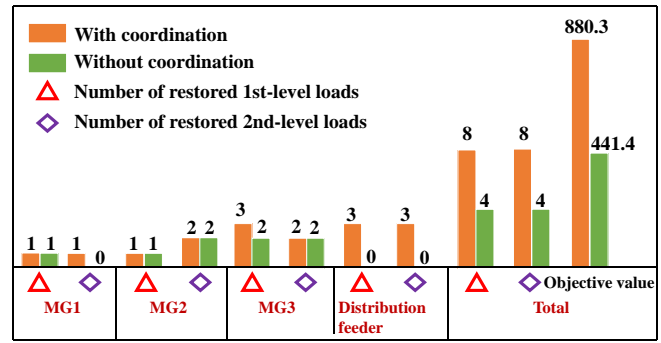

Fig. 7. Comparison between restoration strategies with and without coordinating multiple sources.

Fig. 6. IEEE 123-node Test System with three microgrids.

TABLE IV
DGs AND PVs DATA

| Generator (#) | Microgrid (#) | Node (#) | Generator Type | Rating (kW) |
|---|---|---|---|---|
| DG1 | 1 | 151 | Diesel | 230 |
| DG2 | 1 | 50 | PV | 25 |
| DG3 | 1 | 49 | ES | 35 |
| DG4 | 2 | 300 | Diesel | 140 |
| DG5 | 2 | 108 | Diesel | 70 |
| DG6 | 2 | 105 | PV | 30 |
| DG7 | 2 | 109 | PV | 15 |
| DG8 | 3 | 72 | Diesel | 140 |
| DG9 | 3 | 86 | PV | 170 |

*2) Critical load restoration strategy*

Three MGs are available. Five faults occur and the corresponding faulted lines are 23-25, 35-40, 47-49, 97-98, and 149-150, as shown in Fig. 6.

The proposed method is used to determine the restoration strategy. In the first stage, the MDST is found. By closing all switches except the one on line 54-94, a spanning tree can be formed.

In stage 2, the algorithm converges after three iterations. The rank of all matrices in the solutions are numerically equal to one. The final restoration scheme is shown in Fig. 6. The objective value is 880.3, indicating that eight 1st-level loads and eight 2nd-level loads are restored.

For comparison, the restoration strategy without interconnecting all the MGs is also obtained. Assume that all MGs are in islanded mode and only supply power to their own loads. The comparison is shown in Fig.7. The objective value of strategy by coordinating multiple sources is 99.43% more than that without coordination.

By coordinating multiple sources, the generation capacities of three MGs are fully utilized to restore more critical loads with a higher level of importance. From Fig. 6, it can be seen that MG1 is split into three parts by the faults. The critical load 39 cannot be restored when MG1 operates in the islanded mode. In addition, the generation capacity of MG3 is not sufficient to supply power to all the five critical loads in the islanded mode and the critical load 72 is not restored. MG2 is capable to supply power to all three critical loads within it and has extra power to supply some non-critical loads. However, by coordinating multiple sources, the three MGs can be interconnected with each other through switches and distribution feeders, so the critical loads 39 and 72, and some critical loads on the distribution feeder can be restored by the aggregated generation capacity. Therefore, coordinating multiple sources for service restoration not only contributes to mutual support among MGs but also helps restore critical loads in the distribution system.

*3) Computation time and global optimality*

To test the efficiency and global optimality of the proposed method, a set of 16000 scenarios is randomly generated. Different scenarios have different locations of critical loads and different DGs' power ratings. The simulation results show that the global optimum can be achieved in 99.83% of the scenarios. Since the computation time in mostly consumed by solving semidefinite program in each iteration, the number of iterations is a good indicator of the computation time. The results in Table V show iteration times are between 1 and 3 for all scenarios, indicating that the algorithm can obtain solutions in a reasonable number of iterations. The maximum computation time is 71.56s and the average time is 26.33s. The performance of the method is high and promising for on-line decision support.

TABLE V
THE NUMBER OF SCENARIOS WITH DIFFERENT ITERATION NUMBER.

| Number of Iterations | Number of Scenarios |
|---|---|
| 1 | 14223 |
| 2 | 1722 |
| 3 | 54 |

To compare the proposed algorithm with off-the-shelf optimization solvers for mixed-integer programs, the IEEE 123-node test system is modified to a balanced system. Although there are integer solvers for mixed-integer semidefinite programs, the solvers are still in the testing phase [39], [40]. Therefore, the model for balanced network **CRL-misocp**, given in Appendix A, is used to solve the restoration problem by off-the-shelf optimizers such as the optimizer in MOSEK [37]. For balanced network, the proposed iterative algorithm will solve a SOCP in each iteration. In this paper, 100 scenarios with different locations of critical loads and DGs' power ratings are randomly generated and solved by the proposed algorithm and the mixed-integer optimizer in MOSEK, respectively. The solutions provided by two methods are the same.

*C. Linear Approximation of CLR*

The unbalanced three-phase power flow equation can be modeled by a linear approximation, referred to as LPF [35]. The approximation is realized by assuming that the line losses can be ignored and the voltages are nearly balanced [35]. From the test results of two cases for **CLR**, the two assumptions are satisfied. The restoration model **CLR-milp** using LPF can be formulated and is given in Appendix A. To test the performance of **CLR-milp**, 100 scenarios with different DG capacities and locations of critical loads for unbalanced IEEE 123-node test feeder are generated and solved using the proposed method and **CLR-milp** solvers in MOSEK. The results are in Table VI.

TABLE VI
RESULTS FOR CLR AND CLR-MILP UNDER 100 SCENARIOS.

| Formulations | CLR | CLR-milp |
|---|---|---|
| Average computation time (s) | 23.07 | 1.96 |
| Number of scenarios with same/different results | 22/78 | |

It can be seen that the model **CLR-milp** can be solved more efficiently. However, the results of **CLR-milp** can be over optimistic. In 78 out of the 100 scenarios, one or two additional loads are selected for restoration because the power losses are ignored. The two formulations are compared with more details in our previous work [41].

## VIII. CONCLUSION AND FUTURE WORK

This paper proposes a service restoration method that coordinates multiple sources to restore critical load after extreme events. The critical load restoration problem is solved by a two-stage procedure. The first stage decides the post-restoration topology. In the second stage, the critical load restoration problem is first formulated as a mixed-integer semidefinite program. A relaxed semidefinite program is then proposed by relaxing the unbalanced three-phase power flow constraints and binary variables associated with the load status. An iterative algorithm is proposed to deal with integer variables, which can achieve the global optimality of the primal critical load restoration problem by solving a few semidefinite programs under two conditions. The simulation results validate the effectiveness of the proposed method and indicate that coordinating multiple sources for service restoration can better allocate generation resources to restore more critical loads.

This paper is a start point to explore the benefits of coordinating multiple sources for service restoration and focuses on exploring the first benefit listed in Section III. Additional efforts are required to explore other benefits.

To make full use of limited generation resources, the one snap-shot restoration model **CLR-basic** should be extended to a multi-time-step restoration model to consider operation of electrical islands during the entire outage. Some additional constraints associated with time should be incorporated, such as the limited generation resources (e.g., diesel and gas fuel) constraints [12], the state of charge (SOC) constraints for energy storage systems (ESS), and the ramping rate constraints of DGs [18]. In addition, the scheduling constraints associated with mobile sources can also be considered in the multi-time-step restoration strategy [42]. However, the model will become more complicated because more integer variables will be introduced and the variables in each time step are coupled with those in other steps. Extension of the method proposed in this paper to cope with multi-time-step scenarios is a nontrivial task.

Dynamics and transients induced by switching operations can be significant due to relatively small generation capacity and low inertia of sources. For instance, energizing transformers may result in significant inrush currents [26]. In addition, picking up large loads at once may lead to severe deviation in frequency or even stalling of DGs [43]. The synchronization of MGs with the target island can also be an issue. Therefore, control strategies and transient constraints for optimization should be studied to ensure the succeed in implementing the restorative actions.

Renewables, such as photovoltaic generators (PV), are intermittent source. In order to consider the impacts of renewables on service restoration, an accurate model of uncertainties induced by renewables is needed [44] [45]. In addition, modeling and solution techniques for optimization problems with random variables should be proposed for decision making [20]. Both tasks are challenging and further research is demanded to deal with uncertainties in the critical load restoration problem.

## APPENDIX

*A. Formulations of CLR-misocp and CLR-milp*

The second-order cone model, semidefinite model, and linear approximation model for power flow equations can be used to formulate critical load restoration problem, resulting in three models, i.e., **CLR-misocp** [46], **CLR**, and **CLR-milp** [41], respectively. **CLR** is described in the main body of this paper, while **CLR-misocp** and **CLR-milp** are described here.

*1) CLR-misocp Formulation*

For balanced networks, the critical load restoration problem can be formulated as a mix-integer second-order cone program with consideration of the single-phase power flow [46].

**CLR-misocp:**

$$\max \quad f_{\text{misocp}}(\gamma_i, s_i) = [\sum_{i \in \mathcal{L}} w_i \gamma_i - w_0 \text{real}(\sum_{i \in \mathcal{N}} s_i)] \quad (20)$$

over $\gamma_i \in \{0,1\}$ for $i \in \mathcal{L}$;

$s_i \in \mathbb{C}, v_i \in \mathbb{H}$ for $i \in \mathcal{N}$;

$S_{ij} \in \mathbb{C}, \ell_{ij} \in \mathbb{H}$ for $i \to j$

s.t.

$$\sum_{k:k \to i} S_{ki} - Z_{ki} \ell_{ki} + s_i = \sum_{j:i \to j} S_{ij}, \forall i \in \mathcal{N} \quad (21)$$

$$0 \le s_i + \gamma_i s_{\text{load},i} \le S_{\text{rate},i}, \forall i \in \mathcal{G} \quad (22)$$

$$s_i = -\gamma_i s_{\text{load},i}, \forall i \in \mathcal{N}/\mathcal{G} \tag{23}$$

$$v_{i,\min} \leq v_i \leq v_{i,\max}, \forall i \in \mathcal{N} \tag{24}$$

$$\ell_{ij} \leq |I_{ij,\max}|^2, \forall (i,j) \in \mathcal{E} \tag{25}$$

$$v_j = v_i - \left(S_{ij}Z_{ij}^H + Z_{ij}S_{ij}^H\right) + Z_{ij}\ell_{ij}Z_{ij}^H, \forall i \to j \tag{26}$$

$$\ell_{ij} \geq |S_{ij}|^2/v_i, \forall i \to j \tag{27}$$

where $v_i \coloneqq |V_i|^2$, $\ell_{ij} \coloneqq |I_{ij}|^2$, $S_{ij} \coloneqq V_i I_{ij}^H$, and $Z_{ij}$ is the impedance of line $(i,j)$.

The objective function $f_{\text{misocp}}(\gamma_i, s_i)$ includes the objective of **CLR** and the power losses of the system with a weighting factor $w_0$ to balance the two objectives.

The meanings of constraints (21)-(27) are corresponding to (10)-(16). It can be seen that constraint (27) defines a second-order conic set and it is relaxed from the equation (28), the transformation of the power definition equation [33].

$$\ell_{ij} = |S_{ij}|^2/v_i, \forall i \to j \tag{28}$$

It has been proved that the relaxation is exact when the network is radial and the objective function is convex and the non-decreasing function respect to power injection [34]. The conditions are satisfied in this case, so the power flow results satisfy the constraint (28).

*2) CLR-milp Formulation*

Under two assumptions stated in [35], the SDP-based power flow model can be transformed as a linear approximation of the power flow, referred to as LPF. Based on LPF, the **CLR** problem proposed in this paper can be simplified to a mixed-integer linear program, i.e., **CLR-milp** [41].

**CLR-milp:**

$$\max\ f_{\text{milp}} = f(\boldsymbol{\gamma}) = \sum_{i \in \mathcal{L}} w_i \gamma_i \tag{29}$$

over $\gamma_i \in \{0,1\}$ for $i \in \mathcal{L}$;

$\boldsymbol{s}_i \in \mathbb{C}^{|\alpha_i|}$, $\boldsymbol{V}_i \in \mathbb{H}^{|\alpha_i| \times |\alpha_i|}$ for $i \in \mathcal{N}$;

$\boldsymbol{\Lambda}_{ij} \in \mathbb{C}^{|\alpha_{ij}|}$ for $i \to j$

s.t. (11)- (13)

$$\sum_{k:k \to i} \boldsymbol{\Lambda}_{ki} + \boldsymbol{s}_i = \sum_{j:i \to j} (\boldsymbol{\Lambda}_{ij})^{\alpha_i}, \forall i \in \mathcal{N} \tag{30}$$

$$\boldsymbol{\Lambda}_{ij} \leq \boldsymbol{s}_{ij,\max}, \forall (i,j) \in \mathcal{E} \tag{31}$$

$$\boldsymbol{V}_j = \boldsymbol{V}_i^{\alpha_{ij}} - \left(\delta \cdot \text{DIAG}(\boldsymbol{\Lambda}_{ij})\boldsymbol{Z}_{ij}^H + \boldsymbol{Z}_{ij}\left(\delta \cdot \text{DIAG}(\boldsymbol{\Lambda}_{ij})\right)^H\right), \forall i \to j \tag{32}$$

where $\boldsymbol{\Lambda}_{ij}$ is vector variable comprised of diagonal elements of $\boldsymbol{S}_{ij}$, and $\boldsymbol{S}_{ij}$ can be recovered by $\boldsymbol{S}_{ij} = \delta \cdot \text{DIAG}(\boldsymbol{\Lambda}_{ij})$; function DIAG($\cdot$) returns a diagonal matrix with diagonal elements are the input vector; $\delta = \begin{bmatrix} 1 & \alpha^2 & \alpha \\ \alpha & 1 & \alpha^2 \\ \alpha^2 & \alpha & 1 \end{bmatrix}$ where $\alpha = e^{-\frac{2\pi}{3}j}$.

The objective function of **CLR-milp** is the same as that of **CLR**. Constraints (30)- (32) are simplified from (10), (14) and (15), respectively. The restoration results are the load status and power flow information satisfying all constraints.

*B. Proof of the Optimality Criterion*

To prove *optimality criterion*, it suffices to prove that: 1) The status of loads that are added the constraint $\gamma_i = 1$ is 1 in the global optimal solution; 2) The status of loads those are added the constraint $\gamma_i = 0$ in step 3 is 0 in the global optimal solution; 3) The treatment in step 5 that satisfies the optimality criterion does not influence the global optimality of the primal **CLR**.

For 0th iteration, $W_{\text{CLR-sdp}}^{*[0]}$ is the global optimum after relaxing the integer variable to continuous variable of the primal **CLR**. Obviously, $W_{\text{CLR-sdp}}^{*[0]} \geq W_{\text{CLR}}^*$ and $W_{\gamma_i=1}^{[0]} \leq W_{\text{CLR}}^*$. If $W_{\text{CLR-sdp}}^{*[0]} = W_{\gamma_i=1}^{[0]}$, it indicates that all the load status $\gamma_i^{*[0]}$ in the solution are integers, and this solution is also the global optimal solution of the primal **CLR**. On the other hand, if there is any non-integer load status $\gamma_i^{*[0]}$ in the solution, the relationship will hold that $W_{\text{CLR-sdp}}^{*[0]} > W_{\gamma_i=1}^{[0]}$, indicating there exists at least one non-integer $\gamma_i^{*[0]}$. The latter situation will be discussed in detail.

In the latter situation, $W_{\gamma_i=1}^{[0]}$ can be represented as $W_{\gamma_i=1}^{[0]} = \sum_{k=1}^{K^*} N^k w^k$ ($N^k \in \mathbb{Z}$), where $N^k$ is the number of loads whose $\gamma_i^{*[0]} = 1$ and weighting factor is $w^k$, and $K^*$ satisfies that $w^{K^*} \leq W_{\text{CLR-sdp}}^{*[0]} - W_{\gamma_i=1}^{[0]} < w^{K^*-1}$. Then $W_{\text{CLR-sdp}}^{*[0]}$ can be represented as $W_{\text{CLR-sdp}}^{*[0]} = W_{\gamma_i=1}^{[0]} + \sum_{i:i \in \mathcal{L}_{lk}} (\gamma_i^{*[0]} w^{K^*}) + R$, where $\mathcal{L}_{lk}$ is the set of loads whose weighting factor is $w^{K^*}$ and $\gamma_i^{*[0]}$ is non-integer, and $R$ is the sum of all other terms. In addition, $W_{\text{CLR}}^*$ can be represented as $W_{\text{CLR}}^* = \sum_{k=1}^n N^{k^*} w^k$ ($N^{k^*} \in \mathbb{Z}$).

First, we prove that 1) holds for the first iteration.

The value of objective function is directly affected by the weighting factors and from Section VI. A, we have that for an arbitrary level $j$, $w^j > \sum_{k:k>j} w^k |\mathcal{L}^k|$. Therefore, for the loads whose weighting factor are $w^h$ satisfying $w^h > w^{K^*}$, if there exist loads whose $\gamma_i^{*[0]}$ are non-integers, the cause of the non-integer status is the violation of constraints (13) or (14) rather than limited generation capacities. To clarify this, assume that the load status of one corresponding load in the optimal solution is 1, then it is impossible that the load together with the current fully restored $N^1$ 1st-level loads, $N^2$ 2nd-level loads, ..., and $N^h$ $h$st-level loads are restored fully, i.e., at least one load's status should be 0. To guarantee that the relationship $W_{\text{CLR}}^* > W_{\gamma_i=1}^{[0]}$ holds, at least one load status of loads in all $N^h$ $h$st-level loads should be 0, i.e., the number of loads whose $\gamma_i^{*[0]} = 1$ and weighting factor is $w^h$ should be $n^h$. Because $W_{\text{CLR-sdp}}^{*[0]}$ is the global optimum after relaxing the integer variable, the restoration of $h$-th level and loads in higher level whose $\gamma_i^{*[0]} = 1$ will cause the minimum equivalent power losses to restore more loads in lower levels, in order to optimize the objective. In other words, the effect of $n^h$ loads whose $\gamma_i^{*[0]} = 1$ is the optimal one, so the load status of $h$-th level loads ($w^h > w^{K^*}$) whose $\gamma_i^{*[0]} < 1$ will be 0 in the global optimal solution, and the load status of $h$-th level loads whose $\gamma_i^{*[0]} = 1$ will be 1. Therefore, it is proved that the relationships $n^1 = n^{1^*}, ..., n^{K^*-1} = n^{K^*-1^*}$ hold. For the loads in the $K^*$th level, similarly, the solution $\boldsymbol{\gamma}_{\text{CLR-sdp}}^{*[0]}$ ensures the minimum equivalent power losses of the whole system, so the load status of $K^*$-th level loads whose $\gamma_i^{*[0]} = 1$ will be 1. Therefore, it is proved that 1) holds for the first iteration holds.

Second, we prove 2) holds for the first iteration.

It can be seen that the difference between $W_{\text{CLR-sdp}}^{*[0]}$ and $W_{\gamma_i=1}^{[0]}$ defines an upper bound on the improvement that can be made in the objective value. For the condition that $\mathcal{L}_{c1} \neq \emptyset$, the constraints $\gamma_i = 0$ are added to all the loads in $\mathcal{L}_{c1}$ according to step 3. Consider the situation that there exists an $h$-th level load in the $\mathcal{L}_{c1}$ such that the load status is 1 in the global optimal solution. Then, $W_{\text{CLR}}^* \geq W_{\gamma_i=1}^{[0]} + w^h$ while $W_{\text{CLR-sdp}}^{*[0]} - W_{\gamma_i=1}^{[0]} < w^{K^*-1}$, i.e., $W_{\text{CLR-sdp}}^{*[0]} < W_{\gamma_i=1}^{[0]} + w^h$, which contradicts the known relationship $W_{\text{CLR}}^* \leq W_{\text{CLR-sdp}}^{*[0]}$. Therefore, it is proved that 2) holds for the first iteration.

Finally, we show that 3) holds for the first iteration.

For the condition that step 5 deals with, $n_{\text{re}}$ is the maximum number of loads with $w^{K^*}$ that are estimated to be restored. If the optimality criterion holds, it is indicated that all the $n_{\text{re}}$ loads with $w^{K^*}$ can be fully restored with $\gamma_{\text{CLR},i}^* = 1$. In other words, the criterion indicates that for the primal **CLR**, $(n_{\text{re}} + n^{K^*})$ loads with $w^{K^*}$ can be restored.

*Fact:* Assume that arbitrary combination of $(n_{\text{re}} + N^{K^*})$ loads in all $K^*$-th level loads and loads with higher level whose $\gamma_i^{*[0]} = 1$ can be restored fully, then the maximum difference of power loss between two combinations is so small that cannot be used to restore the minimum loads.

Assume that the global optimum cannot be guaranteed, meaning that at least $(n_{\text{re}} + n^{K^*} + 1)$ loads in $K^*$th level can be fully restored. Then $W_{\text{CLR}}^* \geq W_{\gamma_i=1}^{[0]} + (n_{\text{re}} + 1)w^{K^*}$ while $W_{\text{CLR-sdp}}^{*[0]} - W_{\gamma_i=1}^{[0]} < (n_{\text{re}} + 1)w^{K^*}$, which contradicts with the known relationship that $W_{\text{CLR}}^* \leq W_{\text{CLR-sdp}}^{*[0]}$. Therefore, the treatment in step 5 does not influence the global optimality of the primal **CLR** in the first iteration, i.e., 3) holds for the first iteration.

The sufficient conditions 1)-3) of *optimality criterion* for the latter iterations can be similarly proved.

**Ying Wang** (S'16) received the B.E. degree in electrical engineering from Beijing Jiaotong University, Beijing, China, in 2014. From 2017 to 2018, she was a visiting scholar at Case Western Reserve University, Cleveland, OH, USA. She is currently pursuing her Ph.D. degree at Beijing Jiaotong University, Beijing, China. Her research interests include distribution system restoration and power system resilience.

**Yin Xu** (S'12–M'14–SM'18) received his B.E. and Ph.D. degrees in electrical engineering from Tsinghua University, Beijing, China, in 2008 and 2013, respectively. During 2013-2016, he was an Assistant Research Professor at the School of Electrical Engineering and Computer Science, Washington State University, Pullman, WA, USA. He is currently a Professor at Beijing Jiaotong University, Beijing, China. His research interests include power system resilience, distribution system restoration, and power systems electromagnetic transient simulation.

Dr. Xu is currently serving as Secretary of the Distribution Test Feeder Working Group under the IEEE PES Distribution System Analysis Subcommittee.

**Jinghan He** (M'07–SM'18) received the M.Sc. degree in automation from Tianjin University, Tianjin, China, in 1994 and the Ph.D degree from Beijing Jiaotong University, Beijing, China, in 2007.

She is currently a Professor with Beijing Jiaotong University, Beijing, China. Her research interests include protective relaying, fault distance measurement, and location in power systems.

**Chen-Ching Liu** (S'80–M'83–SM'90–F'94) received the Ph.D. degree from the University of California, Berkeley, CA, USA, in 1983.

He served as a Professor with the University of Washington, Seattle, WA, USA, from 1983–2005. During 2006–2008, he was Palmer Chair Professor with Iowa State University, Ames, IA, USA. During 2008–2011, he was a Professor and Acting/Deputy Principal of the College of Engineering, Mathematical and Physical Sciences with University College Dublin, Ireland. He was Boeing Distinguished Professor of Electrical Engineering and Director of the Energy Systems Innovation Center, Washington State University, Pullman, WA, USA. Currently, he is American Electric Power Professor and Director of the Center for Power and Energy, Virginia Polytechnic Institute and State University, Blacksburg, VA, USA. He is also a Research Professor at Wash-ington State University.

Prof. Liu was the recipient of the IEEE PES Outstanding Power Engineering Educator Award in 2004. In 2013, he received the Doctor Honoris Causa from Polytechnic University of Bucharest, Romania. He served as Chair of the IEEE PES Technical Committee on Power System Analysis, Computing, and Eco-nomics during 2005–2006. He is a Fellow of the IEEE.

**Kevin P. Schneider** (S'00–M'06–SM'08) received his B.S. degree in Physics and M.S. and Ph.D. degrees in Electrical Engineering from the University of Washington, Seattle, WA, USA.

He is currently a Principal Research Engineer with the Pacific Northwest National Laboratory, working at the Battelle Seattle Research Center, Seattle, WA, USA. He is an Adjunct Faculty member with Washington State University, Pullman, WA, USA, and an Affiliate Assistant Professor with the University of Washington, Seattle, WA, USA. His main areas of research are distribution system analysis and power system operations.

Dr. Schneider is a licensed Professional Engineer in Washington State. He is the past Chair of the Distribution System Analysis Sub-Committee and current Secretary of the Power System Analysis, Computing, and Economics (PSACE) Committee.

**Mingguo Hong** (S'93–M'98) received the B.E. degree in electrical engineering from Tsinghua University, Beijing, China; the M.S. degree in mathematics from the University of Minnesota, Duluth, MN, USA; and the Ph. D degree in electrical engineering from the University of Washington, Seattle, WA, USA, in 1991, 1993 and 1998, respectively.

He was a Power System Engineer at ALSTOM Grid, Redmond, WA, USA, from 1998 to 2005; an Associate Professor of Mathematics at Valley Forge Military College, Wayne, PA, USA, from 2005 to 2006; and a Principle Market Engineer a Mid-Continent Independent System Operator, Carmel, IN, USA, from 2006 to 2012. He was an Associate Professor at Case Western Reserve University in Cleveland, OH, USA, from 2012 to 2018. He is currently a Principal Analyst with the ISO New England, located in Holyoke, MA, USA. His research interests are in power system analysis, electricity markets and smart distribution systems.

**Dan T. Ton** received his B.S. degree is in electrical engineering and M.S. degree in business management from the University of Maryland, Baltimore, MD, USA.

He is a Program Manager of Smart Grid R&D within the U.S. Department of Energy (DOE) Office of Electricity Delivery and Energy Reliability (OE). He is responsible for developing and implementing a multi-year R&D program plan for next-generation smart grid technologies to transform the electric grid in the United States through public/private partnerships. Previously, he managed the Renewable Systems Integration program within the DOE Solar Energy Technologies Program.